\documentclass[a4, 12pt]{article}
\usepackage{amsmath, amssymb, amsthm, fullpage}

\usepackage{cleveref}

\theoremstyle{definition}
\newtheorem{dfn}{Definition}[section]
\newtheorem{thm}[dfn]{Theorem}
\newtheorem{fact}[dfn]{Fact}
\newtheorem{lem}[dfn]{Lemma}
\newtheorem{cor}[dfn]{Corollary}

\crefname{dfn}{Definition}{Definitions}
\crefname{thm}{Theorem}{Theorems}
\crefname{fact}{Fact}{Facts}
\crefname{lem}{Lemma}{Lemmas}
\crefname{cor}{Corollary}{Corollaries}

\renewcommand{\subset}{\subseteq}
\newcommand{\power}{\wp}

\begin{document}

\title{On the derived models of self-iterable universes}
\author{Takehiko Gappo \\
Department of Mathematics\\
        Rutgers University, \\
        New Brunswick, NJ, 08854 USA\\\\
Grigor Sargsyan\\
Institute of Mathematics\\
Polish Academy of Sciences, \\
Sopot, Poland}
\date{\today}

\maketitle
\begin{abstract}
We show that if the universe is self-iterable and $\kappa$ is an inaccessible limit of Woodin cardinal then $\textsf{AD}_{\mathbb{R}}+``\Theta$ is regular" holds in the derived model at $\kappa$. The proof is fine-structure free, and only assumes basic knowledge of iteration trees and iteration strategies. Our proof can be viewed as the fine-structure free version of the well-known fact that $\textsf{AD}_{\mathbb{R}}+``\Theta$ is regular" is true in the derived models of hod mice that have inaccessible limit of Woodin cardinals (see for example \cite{Sar09}). However, the proof uses a different set of ideas and is more general. 
\end{abstract}

\section{Introduction}

In recent years, the theory $\Theta_{reg}=_{def}\textsf{AD}_{\mathbb{R}}+``\Theta$ is regular" has become a central theory in the study of determinacy axioms. A lot of work has been done on forcing over models of $\Theta_{reg}$ with a goal of obtaining rich models of $\textsf{ZFC}$. Perhaps the program of obtaining $\textsf{ZFC}$ consequences of determinacy by forcing over models of determninacy begins with the work of Steel and Wesep (see \cite{SteellWesep}) who forcing over a  model of $\Theta_{reg}$ obtained the consistency of saturation of the non-stationary ideal on $\omega_1$. Since then this has been the topic of number of papers including  \cite{square}, \cite{Neeman},  \cite{Woodin83}, \cite{WoodinBook} and many others. Perhaps the most striking result proven this way is due to Woodin who showed how to force $\textsf{ZFC}+\textsf{MM}^{++}(c)$ over a model of $\Theta_{reg}$ (see \cite{WoodinBook}).

The consistency of  $\Theta_{reg}$ was first established by Woodin, who in unpublished work showed that one can obtain a model of $\Theta_{reg}$ from a supercompact cardinal and a proper class of Woodin cardinals. One of the prime goals of \cite{Sar09} was to show that $\Theta_{reg}$ has consistency strength weaker than a Woodin cardinal which is a limit of Woodin cardinals. This drastically reduced the consistency strength of $\Theta_{reg}$ but unfortunately the proof of this fact has not been very accessible to a wide audience as it makes heavy use of fine structure and the theory of hod mice as developed in \cite{Sar09} and later continued in \cite{SteelNorm}. 

In this paper, our goal is to present a proof of $\Theta_{reg}$ assuming self-iterability. Our main theorem is \cref{main_theorem}. Its proof does not use any fine structure theory or any theory of hod mice. It does not, however, remove the theory of hod mice from the mild upper bound calculations of $\Theta_{reg}$. The proof of this fact, as presented in \cite{Sar09}, has two components. The first component is \cite[Theorem 5.22]{Sar09} where it is shown how to get a model of $\Theta_{reg}$ from an existence of a certain hod mouse. The second component is \cite[Chapter 6.5]{Sar09} where it is shown how to obtain the desired hod mouse from a Woodin cardinal that is a limit of Woodin cardinals. Our main theorem presents a fine-structure and hod-mouse free proof of \cite[Theorem 5.22]{Sar09} and is stated in much more general terms. It is still unknown how to obtain, via coarse methods (i.e.\ methods not involving hod mice), a model of $\Theta_{reg}$ from Woodin cardinal that is a limit of Woodin cardinals. 

Recall that $A\subseteq\mathbb{R}$ is $\kappa$-universally Baire if there are trees $(T, S)$ on $\omega\times\kappa$ such that for every poset $\mathbb{P}$ of size $<\kappa$ and for every $V$-generic $g\subseteq \mathbb{P}$, $V[g]\models$``$ p[T]=\mathbb{R}\setminus p[S]$'' and $A=(p[T])^{V}$. We then say that $A$ is universally Baire if for all infinite cardinals $\kappa$, $A$ is $\kappa$-universally Baire.  We also say that $(T, S)$ is a $\kappa$-uB representation of $A$.

We let $\Gamma_\kappa=\{A\subset\mathbb{R}\mid A$ is $\kappa$-uB$\}$ for an infinite cardinal $\kappa$ and $\Gamma_\infty=\bigcap_{\kappa}\Gamma_\kappa$. Given $A\in\Gamma_\kappa$, a poset $\mathbb{P}$ of size $<\kappa$ and a $V$-generic $g\subset\mathbb{P}$, we let $A^g$ be the canonical interpretation of $A$ in $V[g]$. More precisely, letting $(T, S)$ be a $\kappa$-uB representation of $A$, $A^g=(p[T])^{V[g]}$.
It is easy to show that any $\kappa$-uB representation gives the same $A^g$.

Suppose that $\lambda$ is a limit of Woodin cardinals.
Let $\text{Coll}(\omega, <\lambda)$ be the Levy collapse and $g\subset\text{Coll}(\omega, <\lambda)$ be generic.
Then set
\[
\mathbb{R}^{*} = \bigcup_{\alpha<\lambda}\mathbb{R}^{V[g\upharpoonright\alpha]}
\]
and
\[
\mathrm{Hom}^{*} = \{A^{*}\subset\mathbb{R}^{*}_{g} \mid \exists\alpha<\lambda\;(A \in \Gamma^{V[g\upharpoonright\alpha]}_{\lambda})\}
\]
where $A^{*} = \bigcup_{\alpha<\lambda} A^{g\upharpoonright\alpha}$ where $A^{g\upharpoonright\alpha}$ is the canonical interpretation of $A$ in $V[g\upharpoonright\alpha]$.
Then $L(\mathbb{R}^{*}, \mathrm{Hom}^{*})$ is called a derived model (of $V$) at $\lambda$.
This model itself depends on $g$, but its first order theory is independent of $g$ because the model is definable in the homogeneous forcing extension.
When we want to clarify that the derived model is defined in $V[g]$, we will write $L(\mathbb{R}^*, \mathrm{Hom}^*)^{V[g]}$.
The proof of the next theorem can be found in \cite{St09}.

\begin{thm}(Woodin)
Let $\lambda$ be a limit of Woodin cardinals and let $L(\mathbb{R}^{*}, \text{Hom}^{*})$ be a derived model at $\lambda$.
Then
\[
L(\mathbb{R}^{*}, \text{Hom}^{*})\models\textsf{AD}^{+}.
\]
\end{thm}

If we assume that $\lambda$ satisfies a stronger large cardinal assumption then a derived model at $\lambda$ will satisfy a stronger determinacy axiom.
The following result is one of such examples. Unfortunately the proof is unpublished.

\begin{thm}[Woodin]
Assume that there is a proper class of Woodin cardinals that are limits of Woodin cardinals, and suppose $\kappa$ is a supercompact cardinal. 
Let $\lambda>\kappa$ be a limit of Woodin cardinals and let $L(\mathbb{R}^{*}, \text{Hom}^{*})$ be a derived model at $\lambda$.
Then
\[
L(\mathbb{R}^{*}, \text{Hom}^{*})\models\textsf{AD}_{\mathbb{R}}+\Theta\text{ is regular}.
\]
\end{thm}

The assumption on $\lambda$ is far from optimal. Sargsyan showed that $\textsf{AD}_{\mathbb{R}}+$``$\Theta$ is regular'' is relatively consistent to the existence a Woodin limit of Woodin cardinals in \cite{Sar09}, using the theory of hod mice.
Here is our main result.

\begin{thm}\label{main_theorem}
Suppose that $V$ is self-iterable.
Let $\lambda$ be an inaccessible limit of Woodin cardinals and let $L(\mathbb{R}^*, \mathrm{Hom}^*)$ be a derived model at $\lambda$.
Then
\[
L(\mathbb{R}^*, \mathrm{Hom}^*)\models\mathsf{AD}_{\mathbb{R}}+\Theta\text{ is regular}.
\]
\end{thm}

\textbf{Acknowledgments:} The second author's work was supported by the NSF Career Award DMS-1352034 and NSF Award DMS-1954149. The authors would like to thank the referee for a long list of corrections that made the paper more readable. 

\section{Self-iterability}

The notion of self-iterability was defined in \cite{SarTr19}.
Roughly speaking, $V$ is self-iterable if it knows its own unique iteration strategy in any set generic extension.
Ordinary mice don't satisfy this assumption in general, but hod mice at limits of their Woodin cardinals do.
In this section, we summarize the basic facts on self-iterability described in \cite{SarTr19}.

Suppose $P$ is a transitive model of \textsf{ZFC}.
In this article, every iteration tree on $P$ is required to use extenders $E$ in $P$ (and its image under an iteration map) such that the length of $E$ equals the strength of $E$ and is inaccessible in $P$.
Recall that an iteration tree is normal if each extender used in it  is applied to the least possible model and the lengths of the extenders are increasing.
The Unique Branch Hypothesis (\textsf{UBH}, see \cite{IT}) says that every normal iteration tree on $V$ has at most one cofinal wellfounded branch.
The generic Unique Branch Hypothesis (\textsf{gUBH}) is the statement that \textsf{UBH} holds in any set generic extension.

Suppose $P$ is a transitive model of \textsf{ZFC} and $\Psi$ is an iteration strategy for $P$.
If $\xi$ is a $P$-cardinal then $P\vert\xi=H^{P}_{\xi}$ and $\Psi_{P\vert\xi}$ is the fragment of $\Psi$ that acts on normal trees based on $\Psi\vert\xi$.
Also, if $\mathcal{T}$ is an iteration tree on $P$ with last model $Q$, then $\Psi_{\mathcal{T}, Q}$ is the iteration strategy for $Q$ defined by $\Psi_{\mathcal{T}, Q}(\mathcal{U})=\Psi(\mathcal{T}^{\frown}\mathcal{U})$.

Given a strong limit cardinal $\kappa$ and $F\subset\text{Ord}$, we define
\[
W^{\Psi, F}_{\kappa}=(H_{\kappa}, \in, F\cap\kappa, P\vert\kappa, \Psi_{P\vert\kappa}\upharpoonright H_{\kappa}).
\]
Now suppose $X\prec W^{\Psi, F}_{\kappa}$.
We let $M_{X}$ be the transitive collapse of $X$ and $\pi_{X}\colon M_X \to H_{\kappa}$ be the inverse of the transitive collapse.
Basically we denote the preimages of objects in $X$ by using $X$ as a subscript, but we write $P_X$ instead of $(P\vert\kappa)_{X}$.
We also denote the $\pi_{X}$-pullback of $\Psi_{P\vert\kappa}$ by $\Lambda_{X}$.

\begin{dfn}
Let $P$ be a transitive model of \textsf{ZFC} and $\Psi$ be an iteration strategy for $P$.
Then $\Psi$ is a generically universally Baire (guB) strategy if there is a formula $\phi(x)$ in the language of set theory together with three relation symbols and $F\subset\text{Ord}$ such that for every inaccessible cardinal $\kappa$ and for every countable $X \prec W^{\Psi, F}_{\kappa}$, whenever
\begin{itemize}
\item $g \in V$ is $M_{X}$-generic for a poset of size $<\kappa_{X}$ and 
\item $\mathcal{T} \in M_{X}[g]$ is such that for some inaccessible $\eta<\kappa_X$, $\mathcal{T}$ is an iteration tree on $P_{X}\vert\eta$,
\end{itemize}
the following hold:
\begin{itemize}
\item if $\mathrm{lh}(\mathcal{T})$ is a limit ordinal and $\mathcal{T}\in\mathrm{dom}(\Lambda_{X})$ then $\Lambda_{X}(\mathcal{T}) \in M_{X}[g]$,
\item $\mathcal{T}$ is according to $\Lambda_{X}$ if and only if $M_{X}[g]\models\phi[\mathcal{T}]$.
\end{itemize}
We say that $(\phi, F)$ is a generic prescription of $\Psi$.
\end{dfn}

The definition of the guB strategy is motivated by the fact that universally Baireness can be characterized by the existence of club many generically correct hulls (see \cite[Lemma 4.1]{St09}).
A subset $F$ of ordinals in a generic prescription codes such a club, which enables us to refer to every elementary submodel $X$ of $H_{\kappa}$ in the definition of guB strategies.

\begin{dfn}
$V$ is self-iterable if $V\models\mathsf{gUBH}$ and $V$ has a guB iteration strategy.
\end{dfn}

The following fact is the most fundamental result in \cite{SarTr19}.

\begin{fact}\label{guB1}
Assume that $V$ is self-iterable and let $\Psi$ be its guB iteration strategy.
Also suppose that $\delta$ is a Woodin cardinal and $\eta\geq\delta$ is an inaccessible cardinal.
Let $g\subset\mathrm{Coll}(\omega, \eta)$ be $V$-generic.
Then, in $V[g]$, there is an Ord-strategy $\Sigma$ for $V\vert\delta$ such that the following hold:
\begin{enumerate}
\item $\Psi_{V\vert\delta}\subset\Sigma$,
\item Letting $\Phi=\Sigma\upharpoonright\mathrm{HC}^{V[g]}$, $V[g]\models$ ``$\Phi$ is universally Baire,''
\item For all $V[g]$-generic $h$, $\Phi^{h}\upharpoonright V[g]\subset\Sigma$ where $\Phi^{h}$ be the canonical extension of $\Phi$ in $V[g*h]$.
\end{enumerate}
\end{fact}

We describe how to define the above Ord-strategy $\Sigma$.
Let $(\phi, F)$ be a generic prescription of $\Psi$ and take any inaccessible $\lambda>\eta$.
Working in $V\vert\lambda[g]$, we can canonically define a strategy $\Sigma_{\lambda}$ for $V\vert\delta$ that acts on iteration trees in $V\vert\lambda[g]$ as follows:
\begin{itemize}
\item $\mathcal{T}\in\text{dom}(\Sigma_{\lambda})$ if and only if $\mathrm{lh}(\mathcal{T})$ is limit and for every limit $\alpha<\mathrm{lh}(\mathcal{T})$  $V\vert\lambda[g]\models\phi[\mathcal{T}\upharpoonright\alpha^{\frown}\{b\}]\land\forall\beta<\alpha\;\phi[\mathcal{T}\upharpoonright\beta]$, where $b=[0, \alpha)_{\mathcal{T}}$.
\item $\Sigma_{\lambda}(\mathcal{T})=b$ if and only if $V\vert\lambda[g]\models\phi[\mathcal{T}^{\frown}b]\land\forall\alpha<\text{lh}(\mathcal{T})\;\phi[\mathcal{T}\upharpoonright\alpha]$.
\end{itemize}
Then we can show that $\Sigma_{\lambda}$ witnesses the claim holds in $V\vert\lambda[g]$ and that if $\lambda_0<\lambda_1$ are inaccessible cardinals then $\Sigma_{\lambda_0}\subset\Sigma_{\lambda_1}$ by elementary submodel argument.
Then $\Sigma=\bigcup_{\lambda}\Sigma_{\lambda}$ is the desired Ord-strategy.

This strategy $\Sigma$ is called the generic interpretation of $\Psi_{V\vert\delta}$ in $V[g]$ and denoted by $\Psi^{g}_{V\vert\delta}$.
We can also define the generic interpretation of $\Psi_{V\vert\delta}$ in any generic extension as follows:
let $g$ be $V$-generic.
Take an inaccessible cardinal $\eta\geq\delta$ and a $V$-generic $h\subset\mathrm{Coll}(\omega, \eta)$ such that $V[g] \subset V[h]$.
Then set $\Psi^g_{V\vert\delta}=\Psi^h_{V\vert\delta}\upharpoonright V[g]$.
We can show that if $\eta_0\leq\eta_1$ are inaccessible cardinals $\geq\delta$ and $h_0\subset\mathrm{Coll}(\omega, \eta_0)$ and $h_1\subset\mathrm{Coll}(\omega, \eta_1)$ are $V$-generic, then $\Psi^{h_1}_{V\vert\delta}\upharpoonright V[h_0] \subset \Psi^{h_0}_{V\vert\delta}$.
So our definition of $\Psi^g_{V\vert\delta}$ is independent of the choice of $\eta$ and $h$.
Using this notation, the clause (3) of \cref{guB1} can be a little bit strengthened:

\begin{fact}\label{guB2}
Assume that $V$ is self-iterable and let $\Psi, \delta, \eta$ and $g$ as in \cref{guB1}.
Also let $\Phi=\Psi^g_{V\vert\delta}\upharpoonright\text{HC}^{V[g]}$.
Suppose that $\lambda>\eta$ is an inaccessible cardinal and $h$ is $V[g]$-generic.
Let $\Phi^h$ be the canonical interpretation of $\Phi$ in $V[g*h]$.
Then $\Phi^h\subset \Psi^{g*h}_{V\vert\delta}$.
\end{fact}

We need one more fact from \cite{SarTr19}.

\begin{fact}\label{guB3}
Suppose that $V$ is self-iterable and let $\Psi$ be its guB iteration strategy.
Suppose that $\eta<\delta<\lambda$ are ordinals, $\delta$ is a Woodin cardinal and $\lambda$ is an inaccessible cardinal.
To simplify our notation, we write $\Psi$ instead of $\Psi_{V\vert\delta}\upharpoonright V\vert\lambda$.
Suppose that $g$ is $V$-generic for a poset in $V_{\eta}$ and $i\colon V\to P$ is an iteration embedding via $\mathcal{T}\in V[g]$ of length $<\lambda$ based on $V\vert\delta$ according to $\Psi^g$.
Then $i(\Psi)=(\Psi^{g})_{\mathcal{T}, P\vert i(\delta)}\upharpoonright P$.
\end{fact}

\section{Proof of the main theorem}

Now we prove our main theorem \cref{main_theorem}.
The proof makes heavy use of Woodin's extender algebra, so we clarify what is needed here.

\begin{fact}
Let $M$ be a transitive model of \textsf{ZFC} and let $\Sigma$ be an $\omega_1+1$-iteration strategy for $M$.
Suppose that $\kappa<\delta$ are countable ordinals and $M\models$``$\delta$ is a Woodin cardinal.''
Then there is a poset $\mathbb{P}$ which is definable over $M\vert\delta$ such that the following hold:
\begin{itemize}
\item $M\models$``$\mathbb{P}$ has the $\delta$-c.c.''
\item For any $x\in\mathbb{R}$, there is an iteration tree $\mathcal{T}$ on $M$ according to $\Sigma$ based on the window $(\kappa, \delta)$ (i.e., based on $V\vert\delta$ above $\kappa$) with last model $N$ such that $x$ is $(N, \pi^{\mathcal{T}}(\mathbb{P}))$-generic, where $\pi^{\mathcal{T}}\colon M\to N$ is the iteration embedding via $\mathcal{T}$.
\end{itemize}
This $\mathbb{P}$ is called Woodin's extender algebra of $M$ in the window $(\kappa, \delta)$ and denoted by $\mathsf{EA}^{M}_{(\kappa, \delta)}$.
\end{fact}

The interested reader is referred to \cite{Farah}.
We also use the next folklore result proved in \cite{SarTr19}.

\begin{dfn}
$(M, \delta, \Sigma)$ Suslin-co-Suslin captures $A\subset\mathbb{R}$ if $M$ is a countable transitive model of enough large fragment of $\mathsf{ZFC}$, $M\models$``$\delta$ is a Woodin cardinal,'' $\Sigma$ is an $\omega_1$-strategy for $M$, and there are trees $T$ and $S$ in $M$ such that $M\models$``$(T, S)$ is $\delta$-complementing'' and for any $x\in\mathbb{R}$, $x\in A$ if and only if there is an iteration tree $\mathcal{U}$ on $M$ according to $\Sigma$ with last model $N$ such that $x$ is $(N, \mathsf{EA}^N_{\pi^{\mathcal{U}}((\kappa, \delta))})$-generic for some $\kappa<\delta$ and $x\in p[\pi^{\mathcal{U}}(T)]$.
\end{dfn}

\begin{fact}\label{ScS_capturing}
Suppose that $V$ is self-iterable and $\Psi$ is its guB strategy.
Let $\eta<\delta<\lambda$ be ordinals and suppose that $\delta$ is a Woodin cardinal and $\lambda$ is inaccessible.
Let $g$ be $V$-generic for a poset in $V_{\eta}$ and let $A\in\Gamma^{V[g]}_{\lambda}$ witnessed by $(T, S)$.
Then for any countable $X\prec H_{\lambda}[g]$ with $\delta, T, S\in X$, letting $\Lambda_X$ be the $\pi_X$-pullback of the fragment of $\Psi^g$ acting on the iteration trees based on the window $(\eta, \delta)$, $(M_X, \delta_X, \Lambda_X)$ Suslin-co-Suslin captures $A$.
\end{fact}

From now, suppose that $V$ is self-iterable, $\Psi$ is its guB iteration strategy, and $\lambda$ is an inaccessible limit of Woodin cardinals.
Let $\langle\delta_{i} \mid i<\lambda\rangle$ be an increasing enumeration of Woodin cardinals and limits of Woodin cardinals below $\lambda$ and for each $i<\lambda$ and let $w_i$ be the window $(\delta_i, \delta_{i+1})$.
Also, we write $\Psi_{i}$ for the fragment of $\Psi$ acting on iteration trees based on the window $w_i$.

Let $g\subset\mathrm{Coll}(\omega, \lambda)$ be $V$-generic.
We work in $V[g]$, but note that $V\vert\xi$ means $H^{V}_{\xi}$, not $H^{V[g]}_{\xi}$.
Let $\Phi^{g}_{i}=\Psi^{g}_{i}\upharpoonright\text{HC}^{V[g]}$.
Note that this is a universally Baire set in $V[g]$ by \cref{guB1}.
So we can define the canonical interpretation of $\Phi_i$ in the symmetric collapse $V(\mathbb{R}^{*})$ by
\[
\Phi^*_i=\bigcup_{\delta_{i+1}\leq\alpha<\lambda}\Phi^{g\upharpoonright\alpha}_{i}.
\]
By \cref{guB2} (or the argument in the paragraph before that), this is an iteration strategy for $V\vert\delta_i$ that acts on iteration trees in $\mathrm{HC}^*=\bigcup_{\alpha<\lambda}\mathrm{HC}^{V[g\upharpoonright\alpha]}$.
Each $\Phi^*_{i}$ can be canonically coded into a subset of $\mathbb{R}^*$, so we confuse the strategy with its code.
By \cref{guB1}(2), each $\Phi^*_i$ is in $\mathrm{Hom}^*$ for any $i<\lambda$.
Also, it follows from \cref{ScS_capturing} that the Wadge ranks of $\Phi^*_i$'s are cofinal in $\mathrm{Hom}^*$; take any $A^* \in \mathrm{Hom}^*$ and suppose that $A \in \Gamma^{V[g\upharpoonright\eta]}_{\lambda}$ for some $\eta<\lambda$.
Let $i$ be such that $\eta<\delta_i$.
Then \cref{ScS_capturing} gives us the projective definition of $A^{g\upharpoonright\eta}$ in $\Phi^{g\upharpoonright\alpha}_i$ for any $\delta_{i+1}\leq\alpha<\lambda$.
Thus, $A^*$ is projective in $\Phi^*_i$.

Our main lemma, which is not written in \cite{SarTr19}, is the following.

\begin{lem}\label{main_lemma}
For any $i<\gamma$, $\Phi^*_{i+1}$ is not ordinal definable from $\Phi^*_i$ and a real in $L(\mathbb{R}^*, \text{Hom}^*)$.
\end{lem}
\begin{proof}

Suppose toward contradiction that $\Phi^*_{i+1}\in \mathrm{OD}(\Phi^*_i, x)^{L(\mathbb{R}^*, \mathrm{Hom}^*)}$ for some $i<\lambda$ and $x\in\mathbb{R}^*$.
Let $\mathcal{T}$ be the iteration tree according to $\Psi_{i+1}$ with last model $W$ such that $x$ is $(W, \mathsf{EA}^W_{\pi^{\mathcal{T}}(w_{i+1})})$-generic.
By elementarity, $\pi^{\mathcal{T}}(\delta_{i+1})$ is a regular cardinal in $W$.
Since the extender algebra $\mathsf{EA}^W_{\pi^{\mathcal{T}}(w_{i+1})}$ is $\pi^{\mathcal{T}}(\delta_{i+1})$-c.c. in $W$,\\\\
(0) $\pi^{\mathcal{T}}(\delta_{i+1})$ is a regular cardinal in $W[x]$.\\

It follows from the self-iterability that we can recover the derived model of $V$ using Woodin's extender algebra as follows.
Let $\langle a_j \mid j<\lambda'\rangle$ be an enumeration of $\mathbb{R}^*$ in $V[g]$, where $\lambda'$ is the ordinal such that $i+2+\lambda'=\lambda$.
We can inductively construct $\langle W_j, \mathcal{U}_j \mid j<\lambda'\rangle$ such that
\begin{itemize}
\item $W_0=W$.
\item If $j+1<\lambda'$, then $W_j$ is the last model of $\mathcal{U}_{j}$.
\item If $j<\lambda'$ is limit, then $W_j$ is the wellfounded direct limit model along the unique cofinal branch through $\bigoplus_{k<j}\mathcal{U}_{k}$.
\item For $j<\lambda'$, $\mathcal{U}_j$ is an iteration tree on $W_j$ according to $(\Psi^{g}_{i+2+j})_{\mathcal{S}_j, W_j}$, where $\mathcal{S}_{j}=\mathcal{T}^{\frown}\bigoplus_{k<j}\mathcal{U}_{k}$, making $a_j$ $(W_j, \mathsf{EA}^{W_j}_{\pi(w_{i+2+j})})$-generic.
\end{itemize}
Then we set $\mathcal{U}=\bigoplus_{j<\lambda'}\mathcal{U}_j$ and $\mathcal{S}=\mathcal{T}^{\frown}\mathcal{U}$.
Notice that $\mathcal{S}$ is a normal iteration tree on $V$ according to $\Psi^{g}$.
Let $\overline{W}$ be the wellfounded direct model along the unique cofinal branch through $\mathcal{S}$.
Since $\lambda$ is inaccessible, $\pi^{\mathcal{S}}(\lambda)=\lambda$.
By construction, we can pick a $(\overline{W}, \mathrm{Coll}(\omega, <\lambda))$-generic $h$ such that $(\mathbb{R}^*)^{V[g]}=(\mathbb{R}^*)^{\overline{W}[h]}$.
Now we claim that
\[
L(\mathbb{R}^*, \text{Hom}^*)^{V[g]}=L(\mathbb{R}^*, \text{Hom}^*)^{\bar{W}[h]}.
\]
By construction, $(\mathbb{R}^*)^{V[g]}=(\mathbb{R}^*)^{\bar{W}[h]}$.
Also, since\\\\
(1) $\Phi^*_i$'s are cofinal in $(\text{Hom}^*)^{V[g]}$ by \ref{ScS_capturing} and \\
(2) $\pi^{\mathcal{S}}(\Psi_i)$ is a fragment of $(\Psi_i)_{\mathcal{S}, \bar{W}}$ for all $i<\gamma$ by \ref{guB3},\\\\ $(\text{Hom}^*)^{V[g]}=(\text{Hom}^*)^{\bar{W}[h]}$.
Therefore, the two derived models are equal.

Now we have $\Phi^*_{i+1}\in\mathrm{OD}(\Phi^*_i, x)^{L(\mathbb{R}^*, \mathrm{Hom}^*)^{\bar{W}[h]}}\subset\mathrm{OD}(\Phi_i, x)^{\bar{W}[h]}$.
There is $h'\in\bar{W}[h]$ which is $\mathrm{Coll}(\omega, <\lambda)$-generic over $\bar{W}$ such that $\bar{W}[h]=\bar{W}[x][h']$.
So, by homogeneity of $\mathrm{Coll}(\omega, <\lambda)$, $(\Phi^*_{i+1}\upharpoonright \bar{W}[x])\in\bar{W}[x]$ and therefore, $\mathcal{T}\in\bar{W}[x]$. But $\pi^{\mathcal{T}}[\delta_{i+1}]$ is cofinal in $\pi^{\mathcal{T}}(\delta_{i+1})$, and therefore, since $\pi^{\mathcal{T}}[\delta_{i+1}]\in \bar{W}[x]$, $\pi^{\mathcal{T}}(\delta_{i+1})$ is not a regular cardinal in $\bar{W}[x]$. Since the critical point $\mathcal{U}$ is above $\pi^{\mathcal{T}}(\delta_{i+1})$, we get that $\pi^{\mathcal{T}}(\delta_{i+1})$ is not a regular cardinal of $\bar{W}[x]$, contradicting (0).

\end{proof}

If $A$ is a set of reals, $w(A)$ means the Wadge rank of $A$ in the \textsf{AD} context.

\begin{cor}\label{strategies_cofinal_in_theta}
$L(\mathbb{R}^*, \mathrm{Hom}^*)\models\sup_{i<\lambda} w(\Phi^*_i)=\Theta$.
\end{cor}
\begin{proof}
Suppose not.
Take $A\subset\mathbb{R}^*$ such that $A \in L(\mathbb{R}^*, \mathrm{Hom}^*)$ and $w(A)> w(\Phi^*_i)$ for all $i<\lambda$.
Then $A$ is ordinal definable from a $\mathrm{Hom}^*$ set and a real.
(Note that $\mathrm{Hom}^*$ itself is ordinal definable in the derived model because it is equal to all the sets of reals with Wadge rank less than $\alpha$ for some ordinal $\alpha$.) 
Since $\Phi^*_i$'s are cofinal in $\mathrm{Hom}^*$, $A$ is ordinal definable from $\Phi^*_i$ and a real for some $i<\lambda$.
Then it follows that $\Phi^*_{i+1}$ is ordinal definable from $\Phi^*_i$ and a real, which contradicts our main lemma \cref{main_lemma}.
\end{proof}

\begin{cor}\label{AD_R}
$\power(\mathbb{R}^*)\cap L(\mathbb{R}^*, \text{Hom}^*)=\mathrm{Hom}^*$ and $L(\mathbb{R}^*, \mathrm{Hom}^*)\models\mathsf{AD}_{\mathbb{R}}$.
\end{cor}
\begin{proof}
Since $\Phi^*_i$'s are in $\mathrm{Hom}^*$, the first part is an immediate consequence of \cref{strategies_cofinal_in_theta}.
The second part is a direct consequence of Martin and Woodin's result that if \textsf{AD} holds and all sets of reals are Suslin then $\mathsf{AD}_{\mathbb{R}}$ holds.
One can find the proof of this fact in \cite{ADpluesbook}, a manuscript of Larson's book on $\mathsf{AD}^{+}$.
\end{proof}

\begin{cor}\label{DC}
$L(\mathbb{R}^*, \mathrm{Hom}^*)\models\mathsf{DC}$.
\end{cor}
\begin{proof}
Solovay showed in \cite{Solovay} that under $\mathsf{AD}+V=L(\power(\mathbb{R}))$, $\mathsf{DC}$ is equivalent to $\mathrm{cf}(\Theta)>\omega$.
So, it is enough to prove that $L(\mathbb{R}^*, \mathrm{Hom}^*)\models\mathrm{cf}(\Theta)>\omega$.
Suppose toward contradiction that $\langle A_n \mid n<\omega\rangle$ is a sequence of sets of symmetric reals which Wadge ranks are cofinal in $\Theta^{L(\mathbb{R}^*, \mathrm{Hom}^*)}$.
By \cref{AD_R}, each $A_n$ is in $\mathrm{Hom}^*$.
For any $n<\omega$, let $\kappa_n<\lambda$ be such that there is an $B_n \in \mathrm{Hom}^{V[g\upharpoonright\kappa_n]}_{<\lambda}$ such that $A_n=B^*_n$.
Let $\kappa=\sup_{n<\omega}\kappa_n$.
Since $\lambda$ is inaccessible, $\kappa<\lambda$ and there is $i<\lambda$ such that $\delta_i>\kappa$.
By \cref{ScS_capturing}, all $A_n$'s are projective in $\Phi^*_i$.
However, $\Phi^*_{i+1}$ is projective in some $A_n$ and then $\Phi^*_{i+1}$ is projective in $\Phi^*_{i}$, which contradicts \cref{main_lemma}.
\end{proof}

Finally we prove our main result.

\renewcommand{\proofname}{Proof of Theorem 1.3}
\begin{proof}
By \ref{AD_R}, it is enough to show that $\Theta$ is regular in the derived model.
Let $M=L(\mathbb{R}^*, \mathrm{Hom}^*)$.
Note that $\mathbb{R}^M=\mathbb{R}^*$.
Towards a  contradiction suppose that $M\models\Theta$ is singular and let $\Theta=\Theta^{M}$.
Take a cofinal function $f\colon\mathbb{R}^*\to\Theta$ in $M$.
Then $f$ is ordinal definable in $M$ from a set of reals in $\mathrm{Hom}^*$.
Minimizing the ordinal parameter we may assume that $f$ is definable from some $A^*\in\mathrm{Hom}^*$.
Suppose that $\eta<\lambda$ is large enough that $A$ is universally Baire in $V[g\upharpoonright\eta]$ witnessed by trees $T$ and $S$ with $p[T]^{V[g\upharpoonright\eta]}=A^{g\upharpoonright\eta}$.

 Let $\Psi$ be the strategy of $V$. Given $\nu<\lambda$ we let $\Sigma_\nu=\Psi_{V|\nu}$. We let $\Sigma_\nu^*$ be the extension of $\Sigma_\nu$ to $V(\mathbb{R}^*)$.  Let $i_0<\gamma$ be such that $\eta<\delta_{i_0}$ and let $i_1<\gamma$ be such that 
 \begin{center}
     $M\models w(\Sigma^*_{\delta_{i_1}})>\sup f[\mathbb{R}^{V[g\upharpoonright\delta_{i_0}]}]$.
 \end{center}
Such $i_1$ exists because $\mathbb{R}^{V[g\upharpoonright\delta_{i_0}]}$ is countable in $V[g]$ and by \ref{DC} we already know that $\Theta$ has uncountable cofinality in $M$.
Since $f$ is cofinal, there is an $x \in \mathbb{R}^*$ such that\\\\
(1) $M\models f(x)>w(\Sigma^*_{\delta_{i_1}})$.\\

We now find iteration $\mathcal{S}=\mathcal{T}^\frown \mathcal{U}$ as in \ref{main_lemma}, where $\mathcal{T}$ is the iteration to make 
$x$ generic at $\pi^{\mathcal{T}}(\delta_{i_0})$  and $\mathcal{U}$ ``resurrects'' the derived model of $V$. Meaning if $\bar{W}$ is the last model of $\mathcal{S}$ then for some $\bar{W}$-generic $h\subseteq Coll(\omega, <\lambda)$, $M$ is the derived model of $\bar{W}$ at $\lambda$ as computed by $h$ (notice that $\pi^{\mathcal{S}}(\lambda)=\lambda$).
Elementarity of $\pi^{\mathcal{S}}$ implies that\\\\ 
(2) $M\models w((\pi^{\mathcal{S}}(\Sigma_{\delta_{i_1}}))^*) > f[\mathbb{R}^{\bar{W}[h\upharpoonright\pi^{\mathcal{S}}(\delta_{i_0})]}]$.\\\\
Notice that $f$ doesn't change and become something like $\pi^{\mathcal{S}}(f)$. This is because as the critical point of $\pi^{\mathcal{S}}$ is above $\eta$, we can lift $\pi^{\mathcal{S}}$ to $j: V[g\upharpoonright\eta]\rightarrow \bar{W}[g\upharpoonright\eta]$, and we must have that $(j(A))^*=A^*$ (here $j(A)$ makes sense as $A\in V[g\upharpoonright\eta]$). The reason $f$ doesn't change is because $f$ is chosen to be definable in $M$ from $A$. 

However, $(\pi^{\mathcal{S}}(\Sigma_{\delta_{i_1}}))^*$ is a fragment of $\Sigma_{\delta_{i_1}}^*$, so we have that (1) implies \\\\
(3) $M\models f(x)>w((\pi^{\mathcal{S}}(\Sigma_{\delta_{i_1}}))^*)$.\\\\
Clearly, as $x\in\mathbb{R}^{\bar{W}[h\upharpoonright\pi^{\mathcal{S}}(\delta_{i_0})]}$, (3) contradicts (2).

\end{proof}

\bibliographystyle{plain}
\bibliography{final_version}

\end{document}